\documentclass{article}[10pt]
\usepackage{amsmath, amssymb, amsthm, epsf, graphicx, amscd}

\theoremstyle{plain}
\newtheorem{theorem}{Theorem}[section]  
\newtheorem{corollary}[theorem]{Corollary}
\newtheorem{lemma}[theorem]{Lemma}
\theoremstyle{definition}

\theoremstyle{remark}             

\theoremstyle{definition}         

\theoremstyle{definition}         

\numberwithin{equation}{section}

\newcommand{\bb}[1]{\mathbb #1}

\def\c #1{ {\mathcal #1} }

\begin{document}

\title{Equalities, congruences, and quotients of zeta functions \\ in Arithmetic Mirror Symmetry}

\author{C. Douglas Haessig\footnote{University of California, Irvine.}}
\date{March 2006}
\maketitle

\tableofcontents

\section{Introduction}

The purpose of this note is twofold. First, we demonstrate that
under certain conditions we may extend the Arithmetic Mirror Theorem
of \cite[Theorems 1.1 and 6.1]{dW04}. Second, we apply this
extension to the study of the quotient of the zeta functions of
$X_\lambda$ and $Y_\lambda$.

With $\lambda \in \bb C$ we may define a family of complex
projective hypersurfaces $X_\lambda$ in $\bb P^n_{\bb C}$ by
\[
x_1^{n+1} + \cdots + x_{n+1}^{n+1} + \lambda x_1 \cdots x_{n+1} = 0.
\]
With the group
\[
G := \{ (\zeta_1, \ldots, \zeta_{n+1}) | \zeta_i \in \bb C,
\zeta_i^{n+1} = 1, \zeta_1 \cdots \zeta_{n+1} = 1 \}
\]
we may define the (singular) mirror variety $Y_\lambda$ as the
quotient $X_\lambda / G$ where $G$ acts by coordinate
multiplication. It turns out that $Y_\lambda$ is a toric
hypersurface and may be explicitly described as the projective
closure in $\bb P_\Delta$ of the affine toric hypersurface
\[
g(x_1, \ldots, x_n) := x_1 + \cdots + x_n + \frac{1}{x_1 \cdots x_n}
+ \lambda = 0.
\]
Note, $\bb P_\Delta$ is the toric variety obtained from the polytope
in $\bb R^n$ with vertices $\{e_1, \ldots, e_n, -(e_1 + \cdots +
e_n)\}$, where the $e_i$ are the standard basis vectors of $\bb
R^n$. From this description of $Y_\lambda$, if we let $\bb F_q$
denote the finite field with $q$ elements of characteristic $p$, it
makes sense to discuss $\bb F_{q^k}$-rational points of $X_\lambda$
and its mirror $Y_\lambda$ when the parameter $\lambda$ lies in $\bb
F_q$.

When the $gcd(n+1, q^k-1)=1$, there are no $(n+1)$-roots of unity in
the field $\bb F_{q^k}$. Viewing $G$ as a group scheme over $\bb Z$,
this means there are no $\bb F_{q^k}$-rational points of $G$. This
leads one to suspect a direct relation between the $\bb
F_{q^k}$-rational points of $X_\lambda$ and $Y_\lambda$:

\begin{theorem}\label{T: equal}
For every positive integer $k$ such that $gcd(n+1, q^k-1) = 1$, we
have the equality $\# X_\lambda(\bb F_{q^k}) = \#Y_\lambda(\bb
F_{q^k})$.
\end{theorem}

If $W_\lambda$ is a smooth crepant resolution of $Y_\lambda$, then
there is a rational map from $W_\lambda$ to $Y_\lambda$ which is
injective on rational smooth points. Thus, if none of the $\bb
F_{q^k}$-rational points on $Y_\lambda$ are singular points, we see
that $\#Y_\lambda(\bb F_{q^k}) = \#W_\lambda(\bb F_{q^k})$.
Consequently, we have:

\begin{corollary}
Suppose the singular locus of $Y_\lambda$ contains no $\bb
F_{q^k}$-rational points. If $gcd(n+1, q^k-1) = 1$, then we have $\#
X_\lambda(\bb F_{q^k}) = \#Y_\lambda(\bb F_{q^k}) = \# W_\lambda(\bb
F_{q^k})$.
\end{corollary}

Next, when $gcd(n+1, q^k-1) > 1$ we may prove:

\begin{theorem}\label{T: cong}
Let $d := gcd(n+1, q^k-1) > 1$. Then
\begin{enumerate}
\item $\# X_\lambda(\bb F_{q^k}) \equiv 0$ mod $d$,
\item if $n+1$ is a power of a prime $\ell$, then, writing $\lambda
= -(n+1)\psi$ in the new parameter $\psi$, we have
\[
\# X_\lambda(\bb F_{q^k}) \equiv 0 \text{ mod}(\ell d) \quad
\text{and} \quad \# Y_\lambda(\bb F_{q^k}) \equiv
\begin{cases}
1 & \psi^{n+1} = 1 \\
0 & \text{otherwise}
\end{cases}
\text{ mod}(\ell).
\]
\end{enumerate}
\end{theorem}

\noindent Thus, combining \cite[Theorems 1.1]{dW04} and \ref{T:
cong} with the Chinese Remainder Theorem yields:

\begin{corollary}
Suppose $n+1$ is a power of a prime $\ell$ and $gcd(n+1, q) = 1$.
Set $\lambda = -(n+1)\psi$. If $\psi^{n+1} \not= 1$, then for every
positive integer $k$, we have $\# X_\lambda(\bb F_{q^k}) \equiv \#
Y_\lambda(\bb F_{q^k})$ modulo$(\ell q^k)$.
\end{corollary}

Before discussing the proofs of Theorems \ref{T: equal} and \ref{T:
cong}, let us apply Theorem \ref{T: equal} to the quotient of the
zeta functions of $X_\lambda$ and $Y_\lambda$.

\section{Application to zeta functions} From \cite[Theorem
7.3]{dW04}, when $(n+1) | q-1$ then the quotient of the zeta
functions of $X_\lambda$ and $Y_\lambda$, when raised to the
$(-1)^{n}$ power, is a polynomial of specified degree. We suspect
that the divisibility $(n+1)|q-1$ is unnecessary and may be removed
without disturbing the conclusion. Evidence for this is the
following:

\begin{theorem}\label{T: quot}
Let $n+1$ be a prime such that $gcd(n+1,q) =1$. Let $k$ be the
smallest positive integer such that $q^k \equiv 1$ modulo $n+1$.
Assume $X_\lambda$ is non-singular and $\lambda^{n+1} \not=
(-(n+1))^{n+1}$. Then there are positive integers $\rho_1, \ldots,
\rho_s$, each divisible by $k$, and polynomials $Q_1, \ldots, Q_s
\in 1 + T\bb Z[T]$ which are pure of weight $n-3$ and irreducible
over $\bb Z$, such that
\[
\left(\frac{Z(X_\lambda/\bb F_q, T)}{Z(Y_\lambda/\bb F_q, T)}
\right)^{(-1)^n} = Q_1(q^k T^k)^{\rho_1/k} \cdots Q_s(q^k
T^k)^{\rho_s/k}
\]
Furthermore, $\rho_1 + \cdots + \rho_s =
\frac{n(n^n-(-1)^n)}{n+1}-n$. (Note, the polynomials $Q_i$ depend on
$n$ and $\lambda$.)
\end{theorem}

\begin{proof}
For every nonnegative integer $s$ and $j = 1, \ldots, k-1$, we have
$gcd(n+1, q^{sk+j}-1) = 1$. So, by Theorem \ref{T: equal}, we have
$\# X_\lambda(\bb F_{q^{sk+j}}) = \#Y_\lambda(\bb F_{q^{sk+j}})$ for
every $s \geq 0$ and $j = 1, \ldots, k-1$. This implies
\begin{equation}\label{E: quot}
\frac{Z(X_\lambda/\bb F_q, T)}{Z(Y_\lambda/\bb F_q)} = \frac{ \exp
\sum_{s \geq 1} \frac{\# X_\lambda(\bb F_{q^{ks}})}{ks} T^{ks}}{
\exp \sum_{s \geq 1} \frac{\# Y_\lambda(\bb F_{q^{ks}})}{ks} T^{ks}}
= \left( \frac{Z(X_\lambda/\bb F_{q^k}, T^k)}{Z(Y_\lambda/\bb
F_{q^k}, T^k)} \right)^{1/k}
\end{equation}
where the first equality uses the previous sentence and the second
equality is simply definition. By \cite[Theorem 7.3]{dW04}, there
exists a polynomial $R_n(\lambda, T) \in 1 + T\bb Z[T]$ of degree
$\frac{n(n^n-(-1)^n)}{n+1}-n$, pure of weight $n-3$, such that
\[
\left(\frac{Z(X_\lambda/\bb F_{q^k}, T)}{Z(Y_\lambda/\bb F_{q^k},
T)}\right)^{(-1)^n} = R_n(\lambda, q^k T).
\]
Combining this with (\ref{E: quot}) shows us that
\[
\left(\frac{Z(X_\lambda/\bb F_q, T)}{Z(Y_\lambda/\bb F_q,
T)}\right)^{(-1)^n} = R_n(\lambda, q^k T^k)^{1/k}.
\]
Therefore, factorizing $R_n(\lambda, T) = Q_1(T)^{\rho_1} \cdots
Q_s(T)^{\rho_s}$ into irreducibles over $\bb Z$ proves the theorem.
\end{proof}

As a side remark, for $n+1=5$, Theorem \ref{T: quot} explains the
form of the zeta function of $Z(X_\lambda / \bb F_q, T)$ found in
\cite{Ca}.

We note that, for the quintic ($n+1=5$), it follows from
\cite[Equation 10.3]{Ca} and \cite[Equation 10.7]{Ca}, in which they
empirically compute the zeta functions of $X_\lambda$ and
$W_\lambda$, that for $X_\lambda$ smooth,
\[
\frac{Z(X_\lambda/\bb F_q, T)}{Z(W_\lambda/\bb F_q, T)} = R_{\c
A}(q^k T^k, \lambda)^{20/k} R_{\c B}(q^k T^k, \lambda)^{30/k}
\]
where $k$ is the smallest positive integer such that $q^k \equiv 1$
modulo $5$ and the $R$'s are quartic polynomials over $\bb Z$ which
are not necessarily irreducible. Note, $k=1$, 2 or 4. Furthermore,
they have constructed auxiliary curves $\c A$ and $\c B$, both of
genus 4, whose zeta functions experimentally correspond to $R_{\c
A}$ and $R_{\c B}$, respectively. It would be interesting to find
these ``auxiliary varieties'' for general $n+1$ and see how they fit
into the framework of mirror symmetry (if at all).

\section{The proof of Theorem \ref{T: equal}} Without loss, we will
write $q$ instead of $q^k$ in the following proof.

\subsection{Formulas for $X_\lambda(\bb F_q)$ and $Y_\lambda(\bb
F_q)$ in terms of Gauss sums}\label{S: formulas}

Define the Gauss sums $G(k)$ as in Section 2. Also, let $M$ be the
$(n+2)\times(n+2)$-matrix defined in \cite[Section 3]{dW04}.

Define the set
\[
E := \{ k \in \bb Z^{n+2} | 0 \leq k_i \leq q-1 \text{ and } Mk
\equiv 0 \text{ mod}(q-1) \}.
\]
For each $k \in \bb Z^{n+2}$, define $s(k)$ as the number of
non-zero entries in $Mk \in \bb Z^{n+2}$. Next, define
\begin{align*}
E_1 &:= \{ k \in E | \text{not all $k_1, \ldots, k_{n+1}$ are the
same, but }
0 \leq k_{n+2} \leq q-1  \} \\
E_2 &:= \{ k \in E |  k_1 = k_2 = \cdots = k_{n+1}, 0 \leq k_{n+2} \leq q-1 \} \\
E_2^* &:= \{ k \in E_2 | 0 < k_1 < q-1, s(k)=n+2 \} \\
S_k &:= \frac{q^{n+1-s(k)}}{(q-1)^{n+3-s(k)}}\left(
\prod_{j=1}^{n+2} G(k_j) \right) \chi(\lambda)^{k_{n+2}}.
\end{align*}
Now, \cite[Section 3]{dW04} demonstrated that
\[
\#X_\lambda(\bb F_q) = \frac{-1}{q-1} + \sum_{k \in E_1} S_k +
\sum_{k \in E_2} S_k.
\]

Consider $k \in E$. Suppose $k_1 = \cdots = k_{n+1} = 0$. If
$k_{n+2} = 0$, then $S_k = q^{n+1}/(q-1)$, else, if $k_{n+2} = q-1$,
then $S_k = -(q-1)^n$. Similarly, suppose $k_1 = \cdots = k_{n+1} =
q-1$. If $k_{n+2} = 0$, then $S_k = (-1)^{n+1} q^n$, else, if
$k_{n+2} = q-1$, then $S_k = (-1)^n q^{n+1}/(q-1)$.

Next, notice
\begin{equation}\label{E: rows}
Mk = \left(%
\begin{array}{c}
  k_1 + \cdots + k_{n+2} \\
  (n+1)k_1 + k_{n+2} \\
  (n+1)k_2 + k_{n+2} \\
  \vdots \\
  (n+1)k_{n+1} + k_{n+2} \\
\end{array}%
\right) \in  \bb Z^{n+2}.
\end{equation}
If one of the rows equals zero, then we must have $k_i = 0$ for some
$1 \leq i \leq n+1$. Thus, if $k \in E_2$ such that $0 < k_1 < q-1$,
then all the rows of $Mk$ must be non-zero; that is, $s(k) = n+2$.
Putting this together with the last paragraph, we find that for
$\lambda \not= 0$, then
\begin{equation}\label{E: X}
\#X_\lambda(\bb F_q) = \frac{q^{n+1}+(-1)^n q^n - 1 -
(q-1)^{n+1}}{q-1} + \sum_{k \in E_1} S_k + \sum_{k \in E_2^*} S_k.
\end{equation}

If $\lambda = 0$, then Section 2 tells us that $k_{n+2}$ is forced
to equal zero. Thus, in the above calculations, we need to neglect
all terms in which $k_{n+2} \not=0$. Doing this, we obtain
\[
\#X_0(\bb F_q) = \sum_{k \in E_1} S_k + N_0^* + \frac{q^{n+1} -
1}{q-1} + (-1)^{n+1} q^n + \frac{(-1)^n - (q-1)^n}{q}.
\]

Let $N_\lambda^*$ denote the number of $\bb F_q$-rational points on
the affine (toric) variety defined by
\[
g(x_1, \ldots, x_n) := x_1 + \cdots + x_n + \frac{1}{x_1 \cdots x_n}
+ \lambda = 0.
\]
In the proof of \cite[Theorem 5.1]{dW04}, we saw that for $\lambda
\not= 0$
\begin{equation}\label{E: N1}
N_\lambda^*  = \frac{(q-1)^n}{q} + \frac{(-1)^n}{q(q-1)} + \sum_{k
\in E_2^*} S_k.
\end{equation}
Also, if $\lambda = 0$, we may calculate that
\[
N_0^*  = \frac{(q-1)^n}{q} + \frac{(-1)^{n+1}}{q} + \sum_{k \in
E_2^*} S_k.
\]

For ease of reference, recall from \cite[Section 2]{dW04} that, for
all $\lambda \in \bb F_q$), we have
\begin{equation}\label{E: Y}
\#Y_\lambda(\bb F_q) = N_\lambda^* - \frac{(q-1)^n}{q} +
\frac{(-1)^n}{q} + \frac{q^n-1}{q-1}.
\end{equation}

\subsection{Finishing the proof of Theorem \ref{T: equal}}

For $\lambda \not= 0$, combining equations (\ref{E: X}), (\ref{E:
N1}), and (\ref{E: Y}) yields
\[
\#X_\lambda(\bb F_q) - \#Y_\lambda(\bb F_q) = \sum_{k \in E_1} S_k -
(q-1)^n + \frac{q^{n+1} + (-1)^n q^n + (-1)^{n+1} - q^n}{q-1}.
\]
Similarly, if $\lambda = 0$, then
\[
\#X_0(\bb F_q) - \#Y_0(\bb F_q) = q^n[ (-1)^{n+1} + 1] + \sum_{k \in
E_1} S_k.
\]
We may now prove Theorem \ref{T: equal} by demonstrating that the
right-hand sides of the above two formulas equal zero when $gcd(n+1,
q-1) = 1$.

\begin{lemma}
If $gcd(n+1, q-1) = 1$, then the right-hand sides are zero.
\end{lemma}

\begin{proof}
Let $k \in E$. Suppose $k_i \not= 0$ for $1 \leq i \leq q-1$. Then,
by (\ref{E: rows}), we see that $(n+1)k_i \equiv (n+1)k_j$ modulo
$q-1$ for every $1 \leq i,j \leq q-1$. By hypothesis, $n+1$ is
invertible in $\bb Z/(q-1)$, and so $k_i = k_j$. This means that, if
$k \in E_1$ then at least one of the first $n+1$ coordinates must be
zero.

Let $1 \leq i \leq n$. Suppose $k \in E_1$ and its first $i$
coordinates are zero. Then (\ref{E: rows}) tells us that $k_{n+2}$
is either zero or $q-1$. In the \emph{first case}, we have $k_{i+1}
= \cdots = k_{n+1} = q-1$ and $s(k)=n+2$. In the \emph{second case},
again we have $k_{i+1} = \cdots = k_{n+2} = q-1$, but $s(k) =
(n+2)-i$. This leads to the following formulas:
\begin{align*}
&\text{When $k_{n+2}=0$ (first case): } S_k = (-1)^{(n+1)-i} q^n \\
&\text{When $k_{n+2}=q-1$ (second case): } S_k =
(-1)^{n-i}(q-1)^{i-1}q^{(n+1)-i}.
\end{align*}
(Note, if $\lambda = 0$, then $k_{n+2}$ must be zero, and so, the
second case never occurs.) Permuting these zeros around in
$\binom{n+1}{i}$ ways among the first $n+1$ coordinates gives us all
possible points in $E_1$. That is, if we set
\[
A := \sum_{i=1}^n \binom{n+1}{i} (-1)^{(n+1)-i} q^n \qquad
\text{(counts first case)}
\]
and
\[
B := \sum_{i=1}^n \binom{n+1}{i} (-1)^{n-i}(q-1)^{i-1}q^{(n+1)-i}
\qquad \text{(counts second case)},
\]
then we have $\sum_{k \in E_1} S_k = A + B$ for $\lambda \not=0$,
else $\sum_{k \in E_1} S_k = A$ for $\lambda = 0$. Now, by the
binomial theorem, we see that
\[
A = q^n[(1-1)^{n+1} - (-1)^{n+1} - 1] = q^n[(-1)^n - 1]
\]
and
\begin{align*}
B &= (q-1)^{-1}(-1)^n[ (-(q-1) + q)^{n+1} - q^{n+1} -
(-1)^{n+1}(q-1)^{n+1}]\\
&= (q-1)^{-1}[(-1)^n + (-1)^{n+1}q^{n+1} + (q-1)^{n+1}].
\end{align*}
Thus, for $\lambda \not= 0$, we have
\[
\sum_{k \in E_1} S_k = q^n[(-1)^n - 1] + \frac{(-1)^n +
(-1)^{n+1}q^{n+1} + (q-1)^{n+1}}{q-1}
\]
which proves the lemma.
\end{proof}

\section{The proof of Theorem \ref{T: cong}} Let us recall what we
will prove.

\bigskip\noindent{\bf Theorem \ref{T: cong}.} {\it Let $d := gcd(n+1, q^k-1) > 1$. Then
\begin{enumerate}
\item $\# X_\lambda(\bb F_{q^k}) \equiv 0$ modulo $d$.
\item Writing $\lambda = -(n+1)\psi$ in the new parameter $\psi$,
if $n+1$ is a power of a prime $\ell$, then
\[
\# X_\lambda(\bb F_{q^k}) \equiv 0 \text{ mod}(\ell d) \quad
\text{and} \quad \# Y_\lambda(\bb F_{q^k}) \equiv
\begin{cases}
1 & \psi^{n+1} = 1 \\
0 & \text{otherwise}
\end{cases}
\text{ mod}(\ell).
\]
\end{enumerate}}

\begin{proof}
Without loss, we will write $q$ instead of $q^k$ in the following
proof. First, let us prove the congruences on $X_\lambda$. We do
this by gathering all the points of $X_\lambda$ that have the same
number of coordinates zero. For each $1 \leq i \leq n-1$, define
$M_i^*$ as the number of $\bb F_q$-rational points in $\bb
P^{n-i}_{\bb F_q^*}$ which lie on the diagonal hypersurface
\[
x_{i+1}^{n+1} + \cdots + x_{n+1}^{n+1} = 0.
\]
Notice that the group
\[
G_i(\bb F_q) := \{ (\zeta_{i+1}, \ldots, \zeta_{n+1}) | \zeta_j \in
\bb F_q, \zeta_j^{n+1} = 1 \} / \{ (\zeta, \ldots, \zeta) |
\zeta^{n+1} = 1 \}
\]
acts freely on the set of points defining $M_i^*$. Since there are
$d := gcd(n+1, q-1)$ many $(n+1)$-roots of unity in $\bb F_q$, we
have $\# G_i(\bb F_q) = d^{n+1-i} / d = d^{n-i}$. Consequently,
$d^{n-i}$ divides $M_i^*$. Next, let $M_0^*$ be the number of $\bb
F_q$-rational points in $\bb P^n_{\bb F_q^*}$ which lie on
$X_\lambda$. The group
\[
G(\bb F_q) := \{ (\zeta_1, \ldots, \zeta_{n+1}) | \zeta_i \in \bb
F_q, \zeta_i^{n+1} = 1, \zeta_1 \cdots \zeta_{n+1} = 1 \} / \{
(\zeta, \ldots, \zeta) | \zeta^{n+1}=1 \}
\]
acts freely on the points defining $M_0^*$, and so $\# G(\bb F_q) =
d^n$ divides $M_0^*$. Putting this together, we have
\[
\# X_\lambda(\bb F_q) = M_0^* + \sum_{i=1}^{n-1} \binom{n+1}{i}
M_i^*.
\]
This proves the first part of the theorem since each $M_i^*$ is
divisible by $d$. If $n+1$ is a power of a prime $\ell$, then not
only are the $M_i^*$ divisible by $d$, but each of the binomial
factors are divisible by $\ell$; this proves the congruence on
$X_\lambda$ in the second part of the theorem.

We now assume $n+1$ is a power of a prime $\ell$. Let us prove the
congruence on $Y_\lambda$. With $\lambda = -(n+1)\psi$, recall from
(\ref{E: Y}) that
\[
\# Y_\lambda(\bb F_q) = N_\lambda^* - \frac{(q-1)^n}{q} +
\frac{(-1)^n}{q} + \frac{q^n-1}{q-1}
\]
where $N_\lambda^*$ is the number of $\bb F_q$-rational points in
$\bb A^{n+1}_{\bb F_q}$ that satisfy
\begin{equation}\label{E: toric}
\begin{cases}
x_1 + \cdots + x_{n+1} - (n+1)\psi = 0 \\
x_1 \cdots x_{n+1} = 1
\end{cases}.
\end{equation}
We claim that $\# Y_\lambda(\bb F_q) \equiv N_\lambda^*$ modulo
$\ell$. Since $gcd(n+1, q-1)>1$, $q \equiv 1$ modulo $\ell$. Using
the fact that $\frac{q^n-1}{q-1} = q^{n-1} + \cdots + q +1$, we
have: if $\ell$ is an odd prime (the even case is similar), then
\[
- \frac{(q-1)^n}{q} + \frac{(-1)^n}{q} + \frac{q^n-1}{q-1} \equiv -0
+ 1 + n \equiv 0 \quad \text{modulo}(\ell).
\]
This proves the claim.

Since we now have $\# Y_\lambda(\bb F_q) \equiv N_\lambda^*$ modulo
$\ell$, we will concentrate on computing $N_\lambda^*$. Consider
counting the points on (\ref{E: toric}) as follows: suppose a point
$x := (x_1, \ldots, x_{n+1}) \in \bb A^{n+1}_{\bb F_q}$ has two
coordinates equal. Then we may permute these two around in
$\binom{n+1}{2}$ ways without changing the \emph{order} of the other
coordinates. Thus, the orbit of the point $x$ under this type of
permutation contains $\binom{n+1}{2}$ points contained in the affine
toric variety defined by (\ref{E: toric}). Note that we are not
overcounting the points which have multiple pairs of coordinates
being the same, like $(1, 1, 2, 2, 2)$. If all the coordinates of
$x$ are different then we may permute these around in $(n+1)!$ ways.

Putting this together, we find, modulo $\ell$:
\[
N_\lambda^*(\bb F_q) \equiv \#\{x \in \bb A^{n+1}_{\bb F_q} | \text{
all coordinates are equal and $x$ satisfies (\ref{E: toric})} \}.
\]
If all the coordinates are equal then we have the system $(n+1)x -
(n+1)\psi = 0$ and $x^{n+1} = 1$. By hypothesis, $n+1$ is invertible
in $\bb F_q$, thus, we have $x = \psi$ for the first equation, and
so $\psi^{n+1} = 1$ for the second. Therefore,
\[
N_\lambda^* \equiv
\begin{cases}
1 & \psi^{n+1} = 1 \\
0 & \text{otherwise}
\end{cases}
\qquad \text{modulo}(\ell).
\]
\end{proof}

\end{document}